\documentclass{amsart}

\usepackage{amsthm,amsmath,amsfonts,amssymb}
\usepackage{mathrsfs,latexsym,enumitem,enumerate}
\usepackage{graphicx,caption,subcaption}
\usepackage{tkz-euclide}
\newcommand\numberthis{\addtocounter{equation}{1}\tag{\theequation}}

\newtheorem{theorem}{Theorem}

\newtheorem{definition}{Definition}

\begin{document}
\title{The Hausdorff dimension of quasi-circles: a result of Ruelle and Bowen}

\author{Catherine Bruce}
\address{School of Mathematics, University of Manchester, Oxford Road, Manchester M13 9PL, United Kingdom}
\email{catherine-bruce@hotmail.co.uk}

\begin{abstract}
We provide an expanded and clarified proof of the famous result of Bowen and Ruelle giving an asymptotic formula for the Hausdorff dimension of quasi-circles corresponding to the Julia sets of $f(z)=z^2+c$ for small $c$. The proof does not contain new material but has been rewritten to make it more accessible to MSc or PhD students with an interest in dimension theory.
\end{abstract}

\maketitle

\section{Introduction}
\begin{figure}
\includegraphics[width=\textwidth]{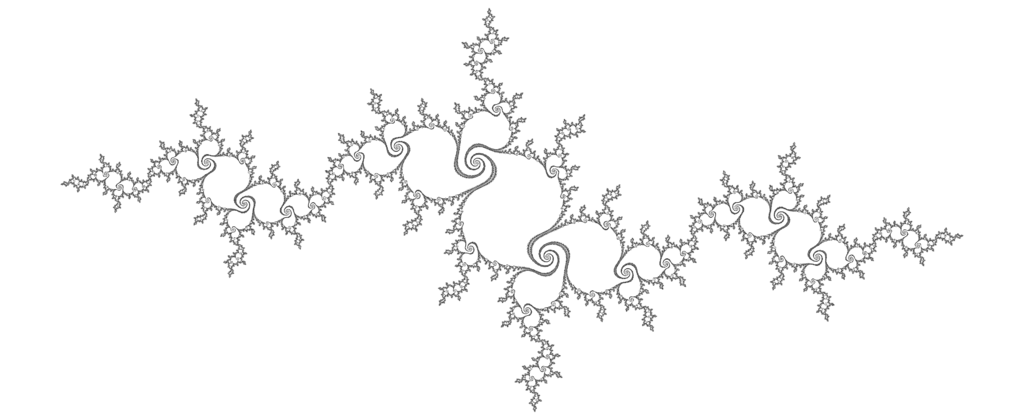}
\caption{Julia set of the quadratic polynomial $f(z)=z^2-1.12+0.222i$}
\label{pic:julia}
\end{figure}

\noindent
Julia sets are an extremely important example of fractals, (see Figure \ref{pic:julia}), which are repellers of certain complex analytic functions. They display approximate self-similarity, as well as being compact, uncountable and perfect. They can be defined in the following way. 
\begin{definition}
The \emph{Julia set} of a function $f$, denoted $J(f)$, is the closure of all repelling periodic points of $f$.
\end{definition}
\noindent For more on the theory of Julia sets including the proofs of the aforementioned properties see \cite{falcfg}, Chapter 14. One of the most fundamental ways to study a fractal set is via its Hausdorff dimension. For a definition and much more on Hausdorff dimension see \cite{falcfg}, Chapter 2.
Julia sets of quadratic functions of the form $f(z)=z^2+c$ have been studied extensively. For $|c|$ large there is an elegant proof, given in \cite{falctifg}, pp228-230, of the following theorem which gives an expression for Hausdorff dimension of such Julia sets. 
\begin{theorem}
Suppose $|c|>\frac{1}{4}(5+2\sqrt{6}).$ Then 
\begin{equation}
\frac{2\log2}{\log4(|c|+|2c|^{\frac{1}{2}})}\le\dim_HJ(f)\le\frac{2\log2}{\log4(|c|-|2c|^{\frac{1}{2}})},
\end{equation}
and so asymptotically we have $\dim_H J(f)\simeq\frac{2\log2}{\log4|c|}.$
\end{theorem} 
The proof views the two branches of the inverse function of $f$ as an iterated function system and uses properties, given in \cite{falcfg}, for finding the Hausdorff dimension of the unique attractor for such a system, which is proved to be the Julia set $J(f)$.

\begin{figure}
\includegraphics[width=5.5cm, height=4.5cm]{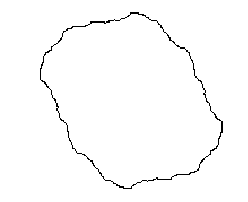}
\caption{A quasi-circle: the Julia set of the quadratic polynomial $f(z)=z^2+\frac{1}{4}i$}
\label{pic:quasi}
\end{figure}

For $|c|$ small, Julia sets of these quadratic polynomials are often referred to as \emph{quasi-circles} (see Figure \ref{pic:quasi}). An asymptotic formula similar to the above is known for small $|c|$, due to Bowen and Ruelle, but the proof is much harder in this case as we are unable to use properties of contractions. Because of this the result is often stated but rarely proved. For example we see the statement in \cite{falcfg}, \cite{mattila}, \cite{mcmullen} and \cite{pesin}. Pollicott also states the result and discusses some detail in \cite{pollicott}, and a detailed proof is given in \cite{zinsmeister}. Our aim is to provide an exposition to ensure Ruelle's extremely important proof is accessible for those studying dimension theory, in particular to MSc and PhD students working in these areas.\\
Let Fix$f^k$ denote the set of all fixed points of $f^k$. Let $S_k\phi(x)=\sum_{j=0}^{k-1}\phi(f^jx)$ be the sum of values of a given potential function $\phi$ at successive iterates of points under $f$. We require a function of $\phi$ of the following form. 
\begin{definition}
We define the \emph{pressure} or \emph{topological pressure} of $\phi$ to be 
\begin{equation}
P(\phi)=\lim_{k\to\infty}\frac{1}{k}\log\sum_{x\in\text{Fix}f^k}\exp(S_k\phi(x)).
\end{equation}
\end{definition}
Pressure, which is one of the most fundamental notions in the \emph{thermodynamic formulism}, is an extremely useful tool for calculating Hausdorff dimension. For more on pressure see \cite{barreira} and \cite{bowentwo}. Let $J(f_c)$ denote a quasi-circle and let $\varphi=-\log|f_c'|$. It can be proved that $P(0)>0$ and $P(s\varphi)\to\infty$ as $s\to\infty$ and thus there is a unique $s$ such that $P(s\varphi)=0$. This and the following theorem are proved by Bowen in \cite{bowen}. This result has been hugely influential and has lead to such results as the analyticity of the Hausdorff dimension as a function of c in the interior of the main cardioid of the Mandelbrot set, originally proved by Ruelle and developed by Smirnov in \cite{smirnov}.  
\begin{theorem}[Bowen's Theorem]
The Hausdorff dimension of a quasi-circle $J(f_c)$ is equal to the unique $s$ such that $P(s\varphi)=0$.
\end{theorem} 
Bowen's result also forms the basis of the proof of our main theorem, which is developed from that given in \cite{ruelle}. 
\begin{theorem}[Ruelle's Theorem]
\label{thm:ruelle}
For $|c|$ small, the Hausdorff dimension of the Julia set $J(f)$ of the map $f:z\to z^2+c$ is equal to
\begin{equation}
\label{eqn:hdcsmall}
1+\frac{|c|^2}{4\log2}+o(|c|^2).
\end{equation}
\end{theorem}
Note that this theorem has been improved on in by Baker and Stallard in \cite{BakerStallard} which gives an estimate of the coefficient in $|c|^3$.
\section{Proof of Theorem \ref{thm:ruelle}}
For $f(z)=z^2+c$ let $\alpha=\frac{1+\sqrt{1-4c}}{2}$ be a fixed point of $f$, i.e. $\alpha^2+c=\alpha$. Note that $\alpha=1-c+O(|c|^2)$ for $|c|$ small. The solution of $f(z)=\eta$ for general $\eta\in \mathbb{C}$ takes the form $\gamma^\epsilon(\eta-c)^{\frac{1}{2}}$ for $\gamma=e^{\pi i}$ and $\epsilon\in \{0,1\}$. Thus we have
\begin{align*}
\{z:f(z)=\alpha\}&=\{\gamma^{\epsilon_1}(\alpha-c)^\frac{1}{2}:\epsilon_1\in\{0,1\}\}\\
\{z:f^2(z)=\alpha\}&=\{\gamma^{\epsilon_2}(\gamma^{\epsilon_1}(\alpha-c)^{\frac{1}{2}}-c)^\frac{1}{2}:\epsilon_2,\epsilon_1\in\{0,1\}\}\\
\vdots\\
\{z:f^n(z)=\alpha\}&=\{\gamma^{\epsilon_n}(...(\gamma^{\epsilon_1}(\alpha-c)^{\frac{1}{2}}-c)^{\frac{1}{2}}...-c)^{\frac{1}{2}}:\epsilon_n,...,\epsilon_1\in\{0,1\}\}.
\end{align*}
Therefore we may encode the set of fixed points $\{z:f^n(z)=\alpha\}$ by dyadic words: for $n\ge 1$ and $\epsilon=\epsilon_1\cdots \epsilon_n \in \{0,1\}^n$ let
\[
\xi(\epsilon)=\gamma^{\epsilon_n}(...(\gamma^{\epsilon_1}(\alpha-c)^{\frac{1}{2}}-c)^{\frac{1}{2}}...-c)^{\frac{1}{2}}.
\]

Fix $s>0$. For $\psi(z)=-s\log|f'(z)|$ we wish to estimate
\begin{align*}
\Delta_n(c,s)&:=\sum_{z:f^n(z)=\alpha}\exp{S_n\psi(z)} \\
&=\sum_{z:f^n(z)=\alpha}\exp\left[{-s\sum_{k=0}^{n-1}\log\left|f'(f^kz)\right|}\right]\\
&=\sum_{\epsilon\in\{0,1\}^n}\exp\left[-s\sum_{k=0}^{n-1}\log2|\xi(\epsilon|_{n-k})|\right],
\end{align*}
where $\epsilon|_{n-k}=\epsilon_1\cdots \epsilon_{n-k}$ denotes the first $n-k$ letters of $\epsilon$, and we have used the fact that  $\left|f'(z)\right|=2|z|$ and $f^k(\xi(\epsilon))=\xi(\epsilon|_{n-k})$.

For $n\ge 1$ and $\epsilon=\epsilon_1\cdots \epsilon_n\in\{0,1\}^n$ we see that
\begin{align*}
\xi(\epsilon)&=\gamma^{\epsilon_n}(...(\gamma^{\epsilon_1}(\alpha-c)^{\frac{1}{2}}-c)^{\frac{1}{2}}...-c)^{\frac{1}{2}}\\
&=\gamma^{\epsilon_n}\left(...\left(\gamma^{\epsilon_2}\gamma^{\frac{\epsilon_1}{2}}\left(\left(\alpha-c\right)^{\frac{1}{2}}-\frac{c}{\gamma^{\epsilon_1}}\right)^{\frac{1}{2}}-c\right)^\frac{1}{2}...-c\right)^{\frac{1}{2}}\\
&=\gamma^{\epsilon_n+\frac{\epsilon_{n-1}}{2}+\cdots+\frac{\epsilon_1}{2^{n-1}}}\left(...\left(\left(\left(\alpha-c\right)^{\frac{1}{2}}-\frac{c}{\gamma^{\epsilon_1}}\right)^{\frac{1}{2}}-\frac{c}{\gamma^{\epsilon_2}\gamma^{\frac{\epsilon_1}{2}}}\right)^\frac{1}{2}...\right.\\
&\qquad\qquad\qquad\qquad\qquad\qquad\qquad\qquad\qquad\qquad\left.-\frac{c}{\gamma^{\epsilon_{n-1}+\frac{\epsilon_{n-2}}{2}+\cdots+\frac{\epsilon_1}{2^{n-2}}}}\right)^{\frac{1}{2}}.\numberthis \label{eqn:xieponeepn}
\end{align*}
Write
$$Q(\epsilon)=\epsilon_n+\frac{\epsilon_{n-1}}{2}+\cdots+\frac{\epsilon_1}{2^{n-1}}$$
and
\begin{equation}\label{eqn2}
r(\epsilon)=\log\left(...\left(\left(\alpha-c\right)^{\frac{1}{2}}-\frac{c}{\gamma^{\epsilon_1}}\right)^{\frac{1}{2}}...-\frac{c}{\gamma^{\epsilon_{n-1}+\cdots+\frac{\epsilon_1}{2^{n-2}}}}\right)^{\frac{1}{2}}
\end{equation}
with the convention $r(\epsilon_1)=\log(\alpha-c)^\frac{1}{2}$ so that
\[
\xi(\epsilon)=\gamma^{Q(\epsilon)}e^{r(\epsilon)}.
\]
From the definition we have
\begin{align*}
r(\epsilon)&=\frac{1}{2}\log\left(e^{r(\epsilon|_{n-1})}-\frac{c}{\gamma^{Q(\epsilon|_{n-1})}}\right)\\
&=\frac{1}{2}\log\left(e^{r(\epsilon|_{n-1})}\left(1-\frac{c}{\xi(\epsilon|_{n-1})}\right)\right)\\
&=\frac{1}{2}r(\epsilon|_{n-1})+\frac{1}{2}\log\left(1-\frac{c}{\xi(\epsilon|_{n-1})}\right).\numberthis\label{eqn:inductiveformforr}
\end{align*}
Now looking at equation (\ref{eqn:xieponeepn}) we see that as $c\to0$ we have $\alpha\to 1$ and 
\begin{equation}\label{eqn1}
\xi(\epsilon|_{n-1})=\gamma^{Q(\epsilon|_{n-1})}+O(n|c|).
\end{equation}
The term $n$ comes from the fact that there are $n-1$ many $(\cdots+\frac{c}{\cdot})^{\frac{1}{2}}$ to develop and each one gives a $\frac{1}{2}$ factor to the first order of $|c|$. By using  $\log(1+z)=z+O(|z|^2)$ and $\frac{1}{1+z}=1+O(|z|)$ when $|z|$ small we have that, to the first order of $|c|$, 
$$r(\epsilon)=\frac{1}{2}r(\epsilon|_{n-1})-\frac{1}{2}\frac{c}{\gamma^{Q(\epsilon|_{n-1})}}+O(n|c|^2).$$
Now, using the above inductive formula for $r(\epsilon|_k)$ up to $k=2$ and use the fact that
\begin{align*}
r(\epsilon_1)=\frac{1}{2}\log \left(\alpha-c\right)&=\frac{1}{2}\log \left(1-2c+O(|c|^2)\right)\\
&=-c+O(|c|^2),
\end{align*}
we eventually see that
$$r(\epsilon)=-c\left(\frac{1}{2\gamma^{Q(\epsilon|_{n-1})}}+\frac{1}{2^2\gamma^{Q(\epsilon|_{n-2})}}+\cdots+\frac{1}{2^{n-1}\gamma^{Q(\epsilon|_1)}}+\frac{1}{2^{n}}\right)+O(n^2|c|^2).$$
Notice
\begin{align*}
\left[\gamma^{{Q(\epsilon)}}\right]^2&=\gamma^{2Q(\epsilon)}\\
&=\exp\left[\left(2\epsilon_n+\epsilon_{n-1}+\cdots+\frac{\epsilon_1}{2^{n-2}}\right)\pi i\right]\\
&=\exp\left[\left(\epsilon_{n-1}+\cdots+\frac{\epsilon_1}{2^{n-2}}\right)\pi i\right]\\
&=\gamma^{Q(\epsilon|_{n-1})}
\end{align*}
since $2\epsilon_n=0\text{ or }2,$ so in either case $e^{2\epsilon_n\pi i}=1.$
Thus
\[
r(\epsilon)=-c\left(\frac{1}{2[\gamma^{Q(\epsilon)}]^2}+\frac{1}{2^2\left[\gamma^{Q(\epsilon)}\right]^{2^2}}+\cdots+\frac{1}{2^{n}\left[\gamma^{Q(\epsilon)}\right]^{2^{n}}}\right)+O(n^2|c|^2).
\]
So, letting $u(\epsilon)=\gamma^{-Q(\epsilon)},$ we have
\begin{align*}
r(\epsilon)=-c\left(\frac{u(\epsilon)^2}{2}+\frac{u(\epsilon)^{2^2}}{2^2}+\cdots+\frac{u(\epsilon)^{2^{n}}}{2^{n}}\right)+O(n^2|c|^2).\numberthis\label{eqn:rassumofpandq}
\end{align*}
Writing $\phi(z)=-\log|f'(z)|=-\log2|z|,$ consider $\phi(\xi(\epsilon)).$ Firstly
\begin{align*}
|\xi(\epsilon)|&=\left|\exp[Q(\epsilon)\cdot\pi i+r(\epsilon)]\right|\\
&=\exp[\operatorname{Re}r(\epsilon)].
\end{align*}
Thus
\begin{align*}
\phi(\xi(\epsilon))&=-\log\left[2\exp\operatorname{Re}r(\epsilon)\right]\\
&=-\log2-\operatorname{Re}r(\epsilon).
\end{align*}
Hence
\begin{align*}
\sum_{k=0}^{n-1}\phi(\xi(\epsilon|_{n-k}))&=-\sum_{k=0}^{n-1}\log2|\xi(\epsilon|_{n-k})|\\
&=-n\log2-\operatorname{Re}\sum_{k=0}^{n-1}r(\epsilon|_{n-k}).\numberthis\label{eqn4}
\end{align*}
Now, summing over $n$ in (\ref{eqn:rassumofpandq}), we see that
\begin{equation}\label{eqn5}
\operatorname{Re}\sum_{k=0}^{n-1}r(\epsilon|_{n-k})=\operatorname{Re}c\Phi_n(u(\epsilon)),
\end{equation}
where
\begin{align*}
\Phi_n(u(\epsilon))&=\frac{1}{2}u(\epsilon)^2+\left(\frac{1+2}{2^2}\right)u(\epsilon)^{2^2}+\left(\frac{1+2+2^2}{2^3}\right)u(\epsilon)^{2^3}+...\\
&\;\;\;\;\;\;\;\;\;\;\;\;\;\;\;\;\;\;\;\;\;\;\;\;\;\;\;\;\;\;\;\;\;\;\;\;\;\;\;\;\qquad\qquad+\left(\frac{1+2+...+2^{n-1}}{2^n}\right)u(\epsilon)^{2^{n}}\\
&=\frac{1}{2}u(\epsilon)^2+\frac{3}{4}u(\epsilon)^4+\frac{7}{8}u(\epsilon)^8+...+\left(1-\frac{1}{2^n}\right)u(\epsilon)^{2^{n}}.
\end{align*}
To second order in $|c|$ we have, using the induction formula (\ref{eqn:inductiveformforr}) \and \eqref{eqn1}, as well as the fact that $r(\epsilon|_{n-1})=O(|c|)$ by \eqref{eqn:rassumofpandq},
\begin{align*}
r(\epsilon)&=\frac{1}{2}r(\epsilon|_{n-1})+\frac{1}{2}\log\left(1-\frac{c}{\xi(\epsilon|_{n-1})}\right)\\
&=\frac{1}{2}r(\epsilon|_{n-1})+\frac{1}{2}\left[-\frac{c}{\xi(\epsilon|_{n-1})}-\frac{1}{2}\left(-\frac{c}{\xi(\epsilon|_{n-1})}\right)^2\right]+O(|c|^3)\\
&=\frac{1}{2}r(\epsilon|_{n-1})-\frac{1}{2}c\exp[-Q(\epsilon|_{n-1})\cdot\pi i-r(\epsilon|_{n-1})]-\frac{c^2}{4(\gamma^{Q(\epsilon|_{n-1})}+O(n|c|))^2}+O(|c|^3)\\
&=\frac{1}{2}r(\epsilon|_{n-1})-\frac{1}{2}cu(\epsilon)^2\exp[-r(\epsilon|_{n-1})]-\frac{1}{4}c^2u(\epsilon)^4+O(n^2|c|^3)\\
&=\frac{1}{2}r(\epsilon|_{n-1})-\frac{1}{2}cu(\epsilon)^2[1-r(\epsilon|_{n-1})]-\frac{1}{4}c^2u(\epsilon)^4+O(n^2|c|^3)\\
&=\frac{1}{2}\left\{r(\epsilon|_{n-1})[1+cu(\epsilon)^2]-cu(\epsilon)^2-\frac{1}{2}c^2u(\epsilon)^4\right\}+O(n^2|c|^3).
\end{align*}
Hence we have
\begin{align*}
\sum_{\epsilon\in\{0,1\}^n}r(\epsilon)&=\frac{1}{2}\sum_{\epsilon\in \{0,1\}^n}\left[r(\epsilon|_{n-1})[1+cu(\epsilon)^2]-cu(\epsilon)^2-\frac{1}{2}c^2u(\epsilon)^4\right]+O(n^2|c|^3)\\
&=\frac{1}{2}\sum_{\epsilon\in \{0,1\}^{n}}r(\epsilon|_{n-1})+O(n^2|c|^3)\\
&=\sum_{\epsilon'\in \{0,1\}^{n-1}}r(\epsilon')+O(n^2|c|^3).\numberthis\label{eqn3}
\end{align*}
where we have used the fact that summing the terms involving $u(\epsilon)^2$ and $u(\epsilon)^4$ over $\epsilon\in\{0,1\}^n$ yields $0$. To see this, take
\[
\sum_{\epsilon\in\{0,1\}^n}r(\epsilon|_{n-1})u(\epsilon)^2=\sum_{\epsilon_n,\epsilon_{n-1}\ldots,\epsilon_1\in\{0,1\}}r(\epsilon_1\ldots \epsilon_{n-1})e^{-(\epsilon_{n-1}-\frac{\epsilon_{n-2}}{2}\cdots-\frac{\epsilon_1}{2^{n-1}})\pi i}
\]
for example. By the definition \eqref{eqn2}, $r(\epsilon_1\ldots \epsilon_{n-1})$ does not depend on $\epsilon_{n-1}$, so summing over $\epsilon_{n-1}\in \{0,1\}$ in the above summation is just summing the same term with different sign, hence $0$. Relation \eqref{eqn3} tells us that,
\begin{align*}
\sum_{\epsilon\in \{0,1\}^n}\sum_{k=0}^{n-1}r(\epsilon|_{n-k})&=\sum_{k=0}^{n-1}\sum_{\epsilon\in \{0,1\}^n}r(\epsilon|_{n-k})\\
&=\sum_{k=0}^{n-1}2^k \sum_{\epsilon'\in \{0,1\}^{n-k}}r(\epsilon')\\
&=\sum_{k=0}^{n-1}2^k \left(O(|c|)+O((n-k)n^2|c|^3)\right)\\
&=O(2^n|c|)+O(2^nn^3|c|^3).\numberthis\label{eqn6}
\end{align*}

Now we move back to the term
\[
\Delta_n(c,s)=\sum_{\epsilon\in\{0,1\}^n}\exp\left[-s\sum_{k=0}^{n-1}\log2|\xi(\epsilon|_{n-k})|\right].
\]
By \eqref{eqn4} we have
\begin{align*}
\Delta_n(c,s)&=\sum_{\epsilon\in\{0,1\}^n}\exp\left\{s\left[-n\log2-\operatorname{Re}\sum_{k=0}^{n-1}r(\epsilon|_{n-k})\right]\right\}\\
&=\sum_{\epsilon\in\{0,1\}^n}2^{-ns}\exp\left[-s\operatorname{Re}\sum_{k=0}^{n-1}r(\epsilon|_{n-k})\right]\\
&=\sum_{\epsilon\in\{0,1\}^n}2^{-ns}\left\{1-s\operatorname{Re}\sum_{k=0}^{n-1}r(\epsilon|_{n-k})+\frac{1}{2}\left[s\operatorname{Re}\sum_{k=0}^{n-1}r(\epsilon|_{n-k})\right]^2 +O(n^3|c|^3)\right\},
\end{align*}
where we have used the fact that $\sum_{k=0}^{n-1}r(\epsilon|_{n-k})=O(n|c|)$ by \eqref{eqn5}. From \eqref{eqn6} we have
\begin{equation}\label{eqn8}
\sum_{\epsilon\in\{0,1\}^n}2^{-ns}\operatorname{Re}\sum_{k=0}^{n-1}r(\epsilon|_{n-k})=O(2^{-n(s-1)}|c|)+O(2^{-n(s-1)}n^3|c|^3).
\end{equation}
For the second last term, by \eqref{eqn5} we have
\begin{align*}
\sum_{\epsilon\in\{0,1\}^n}\left[\operatorname{Re}\sum_{k=0}^{n-1}r(\epsilon|_{n-k})\right]^2 &=\sum_{\epsilon\in\{0,1\}^n}\left[\operatorname{Re}c\Phi_n(u(\epsilon))\right]^2 \\
&= \left[2^{n-1}|c|^2\left(\left(1-\frac{1}{2}\right)^2+\left(1-\frac{1}{2^2}\right)^2+...\right.\right.\\
&\qquad\qquad\qquad\qquad\qquad+\left.\left.\left(1-\frac{1}{2^{n}}\right)^2\right)+O(2^n|c|^2)\right]\\
&=n2^{n-1}|c|^2+o(n)|c|^2+O(2^n|c|^2).\numberthis\label{eqn9}
\end{align*}
Here we have used the following
\begin{equation}\label{eqn7}
\sum_{\epsilon\in\{0,1\}^n}\left(\operatorname{Re}cu(\epsilon)^{2^{r}}\right)\left(\operatorname{Re}cu(\epsilon)^{2^{l}}\right)=0\text{ if }1\le r<l\le n-1
\end{equation}
and
\begin{align*}
\sum_{\epsilon\in\{0,1\}^n}\left(\operatorname{Re}cu(\epsilon)^{2^{r}}\right)^2&=\sum_{\epsilon\in\{0,1\}^n}\frac{1}{2}\left(|c|^2+\operatorname{Re}c^2u(\epsilon)^{2^{r+1}}\right)\numberthis\label{eqn:repusquared}\\
&\left\{
\begin{array}{ccc}
=\frac{1}{2}2^n|c|^2 & \text{if} & r<n-1 \\
=\frac{1}{2}2^n(|c|^2+\operatorname{Re} c^2)& \text{if} & r=n-1\text{ or }n. \\
\end{array}
\right.
\end{align*}
It it the same argument as for \eqref{eqn3} to verify \eqref{eqn7}. To see that (\ref{eqn:repusquared}) holds, consider a complex number $z=x+iy$. Note that $\operatorname{Re}(z^2)=x^2-y^2,$ so
\begin{align*}
(\operatorname{Re}(z))^2&=x^2\\
&=\operatorname{Re}(z^2)+y^2.
\end{align*}
In our case, $z=cu(\epsilon)^{2^{r1}},$ so $|z|=x^2+y^2=|c|$ since $u(\epsilon)$ has modulus 1. So
$$x^2=\operatorname{Re}(z^2)+|c|^2-x^2.$$
Hence
\begin{align*}
(\operatorname{Re}(z))^2&=x^2\\
&=\frac{1}{2}\left(\operatorname{Re}c^2u^{2^{r+1}}+|c|^2\right).
\end{align*}
Gathering \eqref{eqn8} and \eqref{eqn9} we finally obtain
\begin{align*}
\Delta_n(c,s)=2^{-n(s-1)}&+O(2^{-n(s-1)}|c|)+O(2^{-n(s-1)}n^3|c|^3)\\ 
&+\frac{1}{4} s^2 n 2^{-n(s-1)}|c|^2+ o(n)s^22^{-ns}|c|^s+O(2^{-n(s-1)}|c|^2).
\end{align*}
Thus writing $s=1+\beta$ and omitting negligible terms,
$$\Delta_n(c,1+\beta)\approx 2^{-n\beta}\left(1+(1+\beta)^2\frac{|c|^2}{4}n\right),$$
where in this context $\approx$ means that the difference of the two terms is negligible when $|c|\to 0$ and $n\to \infty$. Taking logarithms of the right hand side gives
$$-n\beta\log2+\log\left(1+(1+\beta)^2\frac{|c|^2}{4}n\right)\approx-n\beta\log2+(1+\beta)^2\frac{|c|^2}{4}n.$$
So we have
$$\Delta_n(c,1+\beta)\approx\exp n\left(\frac{|c|^2}{4}(1+\beta)^2-\beta\log2\right).$$
Now take $\beta$ as the positive solution to the equation
\begin{equation}
\label{eqn:solutionofbeta}
\frac{|c|^2}{4}(1+\beta)^2-\beta\log2=0.
\end{equation}
Notice that from \eqref{eqn:solutionofbeta} we have
\[
\beta=\frac{|c|^2}{4 \log 2}(1+\beta)^2=\frac{|c|^2}{4 \log 2}+(2\beta+\beta^2)\frac{|c|^2}{4 \log 2}.
\]
As $\beta\to 0$ as $|c|\to 0$, this gives
$$\beta=\frac{|c|^2}{4\log2}+o(|c|^2).$$
When $\beta$ is the solution of \eqref{eqn:solutionofbeta},
\[
\Delta_n(c,1+\beta)\approx O(1).
\]
Thus the solution of $P(s\phi)=0$ is $1+\beta$ plus lower order terms. By the result of Bowen, we conclude that
\[
\dim_H J(f_c)=1+\frac{|c|^2}{4\log2}+o(|c|^2).
\]
\hfil \qed  

\section*{Acknowledgements}
A special thank you to Xiong Jin without whom this paper would have never come about. Also, to Jonathan Fraser for his invaluable help and support along the way.

\end{document}